\documentclass[a4paper,11pt]{article}
\usepackage{amssymb}
\usepackage{amsthm}
\newtheorem{Lemma}{Lemma}
\newtheorem{Proposition}{Proposition}
\newtheorem{Theorem}{Theorem}

\newcommand{\Aut}{\mathrm{Aut}}

\newcommand{\irr}{\mathrm{irr}}

\newcommand{\Q}{\mathbb{Q}}
\newcommand{\C}{\mathbb{C}}
\newcommand{\bP}{\mathbb{P}}

\newcommand{\sunon}{\overline{S}_{1,n}}
\newcommand{\sgnbar}{\overline{S}_{g,n}}

\title{A remark on the rational cohomology 
of $\overline{S}_{1,n}$}
\author{Gilberto Bini and Claudio Fontanari}
\date{}
\begin{document}
\maketitle

\begin{small}
\begin{center}
\textbf{Abstract}
\end{center}
We focus on the rational cohomology of Cornalba's moduli space 
of spin curves of genus $1$ with $n$ marked points. In particular, 
we show that both its first and its third cohomology group vanish
and the second cohomology group is generated by boundary classes.   
\end{small}

\section{Introduction}
The moduli space of spin curves $\overline{S}_g$ was constructed by 
Cornalba in \cite{Cornalba:89} in order to compactify the 
moduli space of pairs $\{$ smooth genus $g$ complex curve $C$, 
theta-characteristic on $C$ $\}$. 
Cornalba's compactification turns out to be a normal projective variety 
equipped with a finite morphism:
$$
\chi: \overline{S}_g \longrightarrow \overline{\mathcal{M}}_g
$$
onto the Deligne-Mumford moduli space of stable curves of genus $g$
(see \cite{Cornalba:89}, Proposition~(5.2)). 
The geometry of $\overline{S}_g$ (in particular, its Picard group) 
was investigated by Cornalba himself in \cite{Cornalba:89} and in 
\cite{Cornalba:91}; here instead we begin the study of the rational 
cohomology of $\overline{S}_g$. 

As shown by Arbarello and Cornalba in \cite{ArbCor:98}, the rational 
cohomology of $\overline{\mathcal{M}}_g$ vanishes in low odd degree, 
so it seems reasonable to expect that the same holds also for 
$\overline{S}_g$; 
however, a priori it is not clear at all that the morphism $\chi$ does not 
increase cohomology. 
The inductive method of \cite{ArbCor:98} provides indeed an effective tool to 
check our guess, but the set up of the induction requires to 
work with moduli of pointed spin curves. Namely, for all integers 
$g$, $n$ such that $2g-2+n>0$, we consider the moduli spaces
\begin{eqnarray*}
\sgnbar &:=& \{ [(C, p_1, \ldots, p_n; \zeta; \alpha)]: 
(C, p_1, \ldots, p_n) \textrm{ is a genus $g$} \\ 
& & \textrm{quasi-stable projective curve with $n$ marked points}; \\
& &\zeta \textrm{ is a line bundle of degree $g-1$ on $C$ having} \\
& & \textrm{degree $1$ on every exceptional component of $C$, and} \\ 
& &\alpha: \zeta^{\otimes 2} \to \omega_C 
\textrm{ is a homomorphism which} \\ 
& & \textrm{is not zero at a general point of every 
non-exceptional}\\ 
& & \textrm{component of $C$}\}
\end{eqnarray*}

In order to put an analytic structure on $\sgnbar$, we can easily adapt 
Cornalba's construction in \cite{Cornalba:89}: from the universal deformation 
of the stable model of $(C, p_1, \ldots, p_n)$ we obtain exactly as in 
\cite{Cornalba:89}, \S~4, a universal deformation $\mathcal{U}_X \to B_X$ 
of $X = (C, p_1, \ldots, p_n; \zeta; \alpha)$; next, we transplant on 
$\sgnbar$ the structure of $B_X / \Aut(X)$ following \cite{Cornalba:89}, \S~5. 
Alternatively, we can regard $\sgnbar$ as the coarse moduli space associated 
in the easiest case $r=2$ to the stack of $r$-spin curves constructed by 
Jarvis in \cite{Jarvis:00} and revisited by Abramovich and Jarvis in 
\cite{AbrJar:01}. 

We recall that $\sgnbar$ is the union of two connected components, $\sgnbar^+$ and $\sgnbar^-$, which correspond to even and odd theta-characteristics, respectively. The main result of the present paper, which completes the research 
project started in \cite{binfon1} and continued in \cite{binfon2},
is the following:

\begin{Theorem}\label{spinvanishing} 
For every $n$, 
$$
H^1(\sunon^+, \Q) = H^3(\sunon^+, \Q) = 0,
$$
and
$H^2(\sunon^+, \Q)$ is generated by boundary classes.

\end{Theorem}

We note that a similar statement holds true for the moduli space of odd theta-characteristics (see \cite{Getzler:98}) since $\sunon^- \cong {\overline {\mathcal M}}_{1,n}$.
 
In what follows, we work over the field $\C$ of complex numbers; 
all cohomology groups are implicitly assumed to have rational coefficients.

\section{The inductive approach}
As pointed out in the Introduction, we are going to apply the 
inductive strategy developed by Arbarello and Cornalba in 
\cite{ArbCor:98} for the moduli space of curves. 
Namely, we consider the long exact sequence of cohomology with 
compact supports:
\begin{equation}
\label{seq}
\ldots \to H^k_c(S_{1,n}) \to 
H^k(\sunon) \to H^k(\partial S_{1,n})
\to \ldots
\end{equation}
Hence, whenever $H^k_c(S_{1,n}) = 0$,
there is an injection
$H^k(\sunon) \hookrightarrow 
H^k(\partial S_{1,n})$.
Moreover, from \cite{Cornalba:89}, \S~3, it follows that each 
irreducible component of the boundary of $\sunon$ is the image of 
a morphism:
$$
\mu_i: X_i \to \sunon
$$
where either
$$
X_i = {\overline {\mathcal M}}_{0,s+1} \times \overline{S}_{1, t+1}
$$
where $s+t=n$; or
$$
X_i = {\overline{\mathcal M}}_{0,n+2}.
$$

Finally, exactly as in \cite{ArbCor:98}, Lemma~2.6, a bit of Hodge theory 
implies that the map $H^k(\sunon) \to \oplus_i H^k(X_i)$
is injective whenever 
$H^k(\sunon) \to H^k(\partial S_{1,n})$
is. Thus,  we obtain the first claim of Theorem~\ref{spinvanishing} by induction, 
provided we show that 
$H^1_c(S_{1,n}) =
H^3_c(S_{1,n}) = 0$ 
for almost all values of $n$, and we check that 
$H^1(\sunon) = H^3(\sunon) = 0$
for all remaining values of $n$. 
The first task is accomplished by the following
 
\begin{Lemma}\label{harer} We have
$H_k(S_{1,n}) = 0$ for $k > n$.
\end{Lemma}

Indeed, ${\mathcal M}_{1,1} \cong {\mathbb A}^1$ is affine. Moreover, it is well-known that the forgetful morphism ${\mathcal M}_{1,n} \to {\mathcal M}_{1,1}$ is affine. Finally, the morphism $S_{1,n} \to {\mathcal M}_{1,n}$ is finite since it is the restriction of the finite morphism ${\overline {S}_{1,n}} \to {\overline {\mathcal M}}_{1,n}$, hence the claim holds.

Now, we give a closer inspection to $\sunon$. Of course, it is the disjoint union of 
$\overline{S}_{1,n}^+$ 
and 
$\overline{S}_{1,n}^-$,
corresponding to even and odd spin structures respectively. 
However, since the unique odd theta characteristic 
on a smooth elliptic curve $E$ is $\mathcal{O}_E$, 
there is a natural isomorphism
$\overline{S}_{1,n}^- \cong 
\overline{\mathcal{M}}_{1,n}$, so we may restrict our 
attention to $\overline{S}_{1,n}^+$. First of all, the following holds:

\begin{Proposition}
$H^1(\sunon^+)=0$.
\label{accauno}
\end{Proposition}

\proof By the above argument, it is enough to check that $H^1(\sunon^+)$ vanishes for $n=1$. In order to do so, we claim that there is a surjective morphism
$$
f: \overline{\mathcal{M}}_{0,4} \longrightarrow 
\overline{S}_{1,1}^+. 
$$
Indeed, let $(C; p_1, p_2, p_3, p_4)$ be a $4$-pointed stable genus zero 
curve. The morphism $f$ associates to it the admissible covering $E$ 
of $C$ branched at the $p_i$'s, pointed at $q_1$ and equipped with the 
line bundle $\mathcal{O}_E(q_1 - q_2)$, where $q_i$ denotes the point of 
$E$ lying above $p_i$. It follows that 
$$
H^1(\overline{S}_{1,1}^+) \hookrightarrow 
H^1(\overline{\mathcal{M}}_{0,4}) = H^1(\bP^1) = 0   
$$
and Proposition \ref{accauno} is completely proved.

\qed

Recall that the boundary components of $\overline{\mathcal{M}}_{1,n}$ are 
$\Delta_\irr$, whose general member is an irreducible $n$-pointed curve $C$ 
of geometric genus zero with exactly one node, and $\Delta_{1,I}$, whose 
general member is the union of two smooth curves meeting at one node, 
$C_1$ of genus $1$ with marked points labelled by 
$I \subseteq \{ 1, \ldots, n \}$, 
and $C_2$ of genus $0$ with marked points labelled by 
$\{ 1, \ldots, n \} \setminus I$ 
(of course $\vert I \vert \le n-2$). 
The corresponding boundary components of 
$\overline{S}_{1,n}^+$ are:
\begin{itemize}
\item $A^+_\irr$, with an even spin structure on $C$;
\item $B^+_\irr$, with an even spin structure on $C$ blown up at the node;
\item $A_{1,I}^+$, with even theta-characteristics on $C_1$ and $C_2$.
\end{itemize}
Notice that in this case $B_{1,I}^+$, whose general member should carry 
odd theta-characteristics on both $C_1$ and $C_2$, is empty 
since a smooth rational curve has no odd theta-characteristic.

Hence on $\overline{S}_{1,n}^+$
we have the \emph{boundary classes}
$\alpha^+_\irr$, $\beta^+_\irr$, and $\alpha^+_{1,I}$; 
there are also the classes
\begin{eqnarray*}
\delta_\irr &=& p^*(\delta_\irr)\\
\delta_{i,I} &=& p^*(\delta_{i,I})
\end{eqnarray*}
where
$$
p: \overline{S}_{1,n}^+ \to 
\overline{\mathcal{M}}_{1,n}
$$
is the natural projection.
Exactly as in \cite{Cornalba:89}, \S~7, there are relations
\begin{eqnarray}
\label{firstrelation} \delta_\irr &=& \alpha^+_\irr + 2 \beta^+_\irr \\
\label{secondrelation} \delta_{1,I} &=& 2 \alpha^+_{1,I}. 
\end{eqnarray}

\begin{Lemma}\label{h^2spin} 
The vector space $H^2(\overline{S}_{1,2}^+)$
is generated by boundary classes. 
\end{Lemma}

\proof We are going to deduce this from an Euler characteristic computation. 
Indeed, we are going to show that 
\begin{equation}\label{Sbar12}
\chi(\overline{S}_{1,2}^+) = 4.
\end{equation}

Since 
\begin{eqnarray*}
\chi(\overline{S}_{1,2}^+) &=& 
2 h^0(\overline{S}_{1,2}^+) 
- 2 h^1(\overline{S}_{1,2}^+) 
+ h^2(\overline{S}_{1,2}^+) \\
&=& 2 + h^2(\overline{S}_{1,2}^+)
\end{eqnarray*}
from (\ref{Sbar12}) we may deduce that 
$h^2(\overline{S}_{1,2}^+) = 2$.
On the other hand, since the natural projection 
$\overline{S}_{1,2}^+ \to 
\overline{\mathcal{M}}_{1,2}$ is surjective, 
$H^2(\overline{\mathcal{M}}_{1,2})$ injects into 
$H^2(\overline{S}_{1,2}^+)$.
It follows that  
$H^2(\overline{S}_{1,2}^+)$ 
is generated by $\delta_\irr$ and $\delta_{1,\emptyset}$, 
which are linear combinations of $\alpha^+_\irr$, $\beta^+_\irr$,
and $\alpha^+_{1,\emptyset}$ by 
(\ref{firstrelation}) and (\ref{secondrelation}). 

First of all, we compute $\chi(S_{1,1}^+)$. 
It is clear that 
$$
\chi(\overline{S}_{1,1}^+) 
= 2 h^0(\overline{S}_{1,1}^+) 
- h^1(\overline{S}_{1,1}^+) = 2.
$$
On the other hand, $\partial S_{1,1}^+$ 
consists of exactly two points, 
corresponding to a $3$-pointed rational curve with two marked points 
either identified or joined by an exceptional component. Hence 
\begin{equation}\label{S11}
\chi(S_{1,1}^+) = 
\chi(\overline{S}_{1,1}^+) - 
\chi(\partial S_{1,1}^+) = 0.
\end{equation}
Next, we compute $\chi(S_{1,2}^+)$. 
The natural projection $S_{1,2}^+  
\to \mathcal{M}_{1,2}$ is generically three-to-one, 
but there are a few special fibers with less than three points. 
Indeed, let $(E; p_1, p_2)$ be a smooth $2$-pointed elliptic curve.
 
The linear series $\vert 2 p_1 \vert$ provides a realization of $E$ 
as a two-sheeted covering of $\bP^1$ ramified over $\infty$, $0$, $1$ 
and $\lambda$. Denote by $q_0$, $q_1$, and $q_\lambda$ the points of $E$ 
lying above $0$, $1$, and $\lambda$, so that the three even 
theta-characteristics of $E$ are given by $\mathcal{O}_E(p_1 - q_0)$, 
$\mathcal{O}_E(p_1 - q_1)$, and $\mathcal{O}_E(p_1 - q_\lambda)$. 
If $\lambda = \frac{1}{2}$ and $p_2 = q_\lambda$, then 
the projectivity of $\bP^1$ defined by $z \mapsto 1 - z$ induces 
an automorphism of $(E; p_1, p_2)$ exchanging $\mathcal{O}_E(p_1 - q_0)$ 
and $\mathcal{O}_E(p_1 - q_1)$. 
If $\lambda = - \omega$ (with $\omega^3 = 1$) and 
$p_2$ is one point lying above $\frac{\omega}{\omega - 1}$, 
then the projectivity of $\bP^1$ defined by
$z \mapsto \frac{z+\omega}{\omega}$ induces an
automorphism of $(E; p_1, p_2)$ that exchanges ciclically its three 
even theta-characteristics. 
Since it is clear (for instance, from \cite{Hartshorne:77}, 
IV, proof of Corollary~4.7)
that the above ones are the only exceptional cases, we have:
\begin{equation}\label{S12}
\chi(S_{1,2}^+) = 
3 \chi(\mathcal{M}_{1,2} \setminus \{2 \textrm{ points}\})
+ 2 \chi(\textrm{point}) + \chi(\textrm{point}) = 0.  
\end{equation}

In fact, $\chi(\mathcal{M}_{1,2}) = 1$, as observed in 
\cite{ArbCor:98}, (5.4).
Finally we turn to the Euler characteristic of 
$\overline{S}_{1,2}^+$. 
From \cite{Cornalba:89}, Examples (3.2) and (3.3), and 
\cite{ArbCor:98}, Figure~1, we may deduce that
$$
\chi(\overline{S}_{1,2}^+) = 
\chi(S_{1,2}^+) + 
2 \chi(\mathcal{M}'_{0,4}) +
\chi(S_{1,1}^+) + 4,
$$
where $\mathcal{M}'_{0,n}$ denotes the quotient 
of $\mathcal{M}_{0,n}$ modulo the operation of 
interchanging the labelling of two of the marked 
points.

Since $\chi(\mathcal{M}'_{0,4})=0$ (see \cite{ArbCor:98}, (5.4)), 
relation (\ref{Sbar12}) follows from (\ref{S11}) and (\ref{S12}).  

\qed

Let $P$ a finite set with $\vert P \vert = n$ and let $x$ and $y$ be distinct 
and not belonging to $P$; define
$$
\xi: \overline{\mathcal{M}}_{0, P \cup \{ x,y \}} \longrightarrow 
B_\irr^+ \hookrightarrow 
\overline{S}_{1,n}^+
$$
by joining the points labelled $x$ and $y$ with an exceptional component  
and taking the unique even theta characteristic on the resulting curve. 
Then the analogue of Lemma~4.5 in \cite{ArbCor:98} holds:

\begin{Lemma}\label{kernel} 
The kernel of 
$$
\xi^*: H^2(\overline{S}_{1,n}^+)
\longrightarrow  H^2(\overline{\mathcal{M}}_{0, P \cup \{ x,y \}})
$$ 
is one-dimensional and generated by $\delta_\irr$. 
\end{Lemma}

\proof By \cite{ArbCor:98}, Lemma~3.16, it is clear that 
$\xi^*(\delta_\irr)=0$. Moreover, from Lemma~\ref{h^2spin} it follows that 
$H^2(\overline{S}_{1,2}^+)$ is generated by 
$\delta_\irr$ and $\delta_{1,\emptyset}$; since  
$\delta_{1, \emptyset}$ pulls back to $\delta_{ 0, \{ x,y \}}$, 
which is not zero, the claim holds for $n=2$.
Hence we can apply the inductive argument of \cite{ArbCor:98}, pp.~113--114. 
It follows that if $\xi^*(\alpha) = 0$ for $\alpha \in 
H^2(\overline{S}_{1,n}^+)$ 
then there exists a constant $a$ such that $\alpha - a \delta_\irr$ 
restricts to zero on all boundary components of 
$\overline{S}_{1,n}^+$ 
different from $A_\irr^+$. 
However, we claim that $A_\irr^+$ is linearly equivalent to $2 B_\irr^+$. 
Indeed, this is clear in $\overline{S}_{1,1}^+
\cong \bP^1$. If $\pi: \overline{S}_{1,n}^+
\to \overline{S}_{1,1}^+$ is the natural 
forgetful map, then $A_\irr^+ = \pi^*(A_\irr^+)$ and 
$B_\irr^+ = \pi^*(B_\irr^+)$.
Hence $\alpha - a \delta_\irr$ restricts to zero on all boundary components 
of $\overline{S}_{1,n}^+$ and the claim 
follows exactly as in \cite{ArbCor:98}, Lemma~4.5, from 
Proposition~\ref{harer} and the analogue of \cite{ArbCor:98}, Proposition~2.8.

\qed

\begin{Proposition}\label{boundary}
The vector space $H^2(\overline{S}_{1,n}^+)$ 
is generated by boundary classes.
\end{Proposition}

\proof Let $V$ be the subspace of 
$H^2(\overline{S}_{1,n}^+)$
generated by the elements $\alpha_{1,I}^+$. 
In view of Lemma~\ref{kernel} and (\ref{firstrelation}), 
it will be sufficient to show that the morphism $\xi^*$ 
vanishes modulo $V$. The proof is by induction on $n$:
for the inductive step we refer to \cite{ArbCor:98}, pp.~114--118,
while the basis of the induction is provided by Lemma~\ref{h^2spin}.

\qed 

Finally, we are also able to prove the last part of Theorem (\ref{spinvanishing}). 

\begin{Lemma}\label{h^3}
We have $H^3(\overline{S}_{1,n}^+)=0$.
\end{Lemma}

\proof Once again, by the exact sequence (\ref{seq}) and Lemma 1, it is enough to check that $H^3({\overline S}_{1,n}^+)=0$ for $ n \leq 3$. First, we deal with the case $n=3$. By Proposition~\ref{boundary}, 
$H^2(\overline{S}_{1,3}^+)$ 
is generated by the six boundary classes $\alpha_\irr^+$, 
$\beta_\irr^+$, $\alpha_{1, \emptyset}^+$, 
$\alpha_{1, \{ 1 \} }^+$, $\alpha_{1, \{ 2 \} }^+$, and 
$\alpha_{1, \{ 3 \} }^+$. Notice further that 
$\alpha_\irr^+$ and $\beta_\irr^+$ are linearly dependent. 
Indeed, if $\pi: \overline{S}_{1,3}^+
\to \overline{S}_{1,1}^+$ is the natural 
forgetful map, from \cite{ArbCor:98}, Lemma~3.1~(iii) it follows  
that $\alpha_\irr^+ = \pi^*(\alpha_\irr^+)$ and 
$\beta_\irr^+ = \pi^*(\beta_\irr^+)$, while Poincar\'e duality
yields $H^2(\overline{S}_{1,1}^+) \cong 
H^0(\overline{S}_{1,1}^+) \cong \Q$.
Hence we deduce $h^2(\overline{S}_{1,3}^+) \le 5$; 
next, we claim that
\begin{equation}\label{Sbar13}
\chi(\overline{S}_{1,3}^+) = 12.
\end{equation}
The statement is a direct consequence of the claim, since
\begin{eqnarray*}
\chi(\overline{S}_{1,3}^+) &=& 
2 h^0(\overline{S}_{1,3}^+) 
+ 2 h^2(\overline{S}_{1,3}^+) 
- 2 h^1(\overline{S}_{1,3}^+) \\
& & - h^3(\overline{S}_{1,3}^+) 
\le 12 - h^3(\overline{S}_{1,3}^+).
\end{eqnarray*}

First of all, we compute $\chi(S_{1,3}^+)$. 
The natural projection $S_{1,3}^+  
\to \mathcal{M}_{1,3}$ is generically three-to-one, 
but there is a special fiber with only one point. 
Indeed, if $(E; p_1, p_2, p_3)$ is a smooth $3$-pointed elliptic curve
realized by the linear series $\vert 2 p_1 \vert$ 
as a two-sheeted covering of $\bP^1$ ramified over $\infty$, $0$, $1$ 
and $- \omega$ (with $\omega^3 = 1$) and 
$p_2$, $p_3$ are the two points lying above $\frac{\omega}{\omega - 1}$, 
then the projectivity of $\bP^1$ defined by
$z \mapsto \frac{z+\omega}{\omega}$ induces 
automorphisms of $(E; p_1, p_2, p_3)$ exchanging ciclically its three 
even theta-characteristics.
Therefore we have: 
\begin{equation}\label{S13}
\chi(S_{1,3}^+) = 
3 \chi(\mathcal{M}_{1,3} \setminus \{\textrm{ point}\})
+ \chi(\textrm{point}) = - 2.  
\end{equation}

Recall that $\chi(\mathcal{M}_{1,3}) = 0$, as observed in 
\cite{ArbCor:98}, (5.4).
Finally we turn to the Euler characteristic of 
$\overline{S}_{1,3}^+$. 
From \cite{Cornalba:89}, Examples (3.2) and (3.3), and 
\cite{ArbCor:98}, Figure~2, it is clear that
\begin{eqnarray*}
\chi(\overline{S}_{1,3}^+) &=& 
\chi(S_{1,3}^+) + 
2 \chi(\mathcal{M}'_{0,5}) +
\chi(S_{1,1}^+)\chi(\mathcal{M}_{0,4}) + \\
& &
3 \chi(S_{1,2}^+) +
2 \chi(\mathcal{M}_{0,4}) + 
12 \chi(\mathcal{M}'_{0,4}) + 
3 \chi(S_{1,1}^{\hspace{0.05cm}(0),+}) + 14. 
\end{eqnarray*}
Since $\chi(\mathcal{M}_{0,4})=-1$, $\chi(\mathcal{M}'_{0,4})=0$
and $\chi(\mathcal{M}'_{0,5})=1$ (see \cite{ArbCor:98}, (5.4)), 
now (\ref{Sbar13}) follows from (\ref{S11}), (\ref{S12}) and 
(\ref{S13}). Finally, by Hodge theory of complex projective orbifolds, the surjective morphism ${\overline S}_{1,3} \to {\overline S}_{1,n}$ for $n \leq 3$ induces an injective morphism $H^3({\overline S}_{1,n}) \to H^3({\overline S}_{1,3})$ for $n \leq 3$. Hence the claim follows.

\qed

{\bf Acknowledgements.} The first named author is supported by the University 
of Milan ``FIRST'' and by MiUR Cofin 2006 - Variet\`{a} algebriche: geometria, aritmetica, strutture di Hodge. The second named author is supported by 
MiUR Cofin 2006 - Geometria delle variet\`{a} algebriche e dei loro spazi di 
moduli.

\vspace{1cm}

\noindent
Gilberto Bini\\
Dipartimento di Matematica \\ 
Universit\`a degli Studi di Milano \\ 
Via C. Saldini 50 \\
20133 Milano (Italy)\\
e-mail: gilberto.bini@mat.unimi.it

\vspace{0.5cm}

\noindent
Claudio Fontanari\\
Dipartimento di Matematica \\ 
Politecnico di Torino \\
Corso Duca degli Abruzzi 24 \\ 
10129 Torino (Italy)\\
e-mail: claudio.fontanari@polito.it

\end{document}